\RequirePackage{fix-cm}
\documentclass[smallextended]{svjour3}       
\smartqed  
\usepackage{graphicx}
\usepackage{amssymb}
\usepackage{amsmath}
\usepackage{enumitem}

\newcommand{\ra}{\rightarrow}
\newcommand{\N}{\mathbb{N}}
\newcommand{\R}{\mathbb{R}}

\begin{document}

\title{Functions that preserve totally bounded sets vis-\'a-vis stronger notions of continuity
}


\author{Lipsy Gupta         \and
       S. Kundu 
}


\institute{Lipsy Gupta  \at
	Department of
	Mathematics, Indian Institute of Technology Delhi, \\New Delhi, 110016, India\\
	\email{lipsy1247@gmail.com}        
	\and
	S. Kundu \at
	Department of
	Mathematics, Indian Institute of Technology Delhi, \\New Delhi, 110016, India\\
	\email{skundu@maths.iitd.ac.in} 
}

\date{Received: date / Accepted: date}

\maketitle

\begin{abstract}
	A function between two metric spaces is said to be totally bounded regular if it preserves totally bounded sets. These functions need not be continuous in general. Hence the purpose of this article is to study such functions vis-\'a-vis continuous functions and functions that are stronger than the continuous functions such as Cauchy continuous functions, some Lipschitz-type functions etc. We also present some analysis on strongly uniformly continuous functions which were first introduced in \cite{[BL2]} and study when these functions are stable under reciprocation. 
\end{abstract}
\section{Introduction}
A subset $A$ of a metric space $(X,d)$ is called totally bounded if for every $\epsilon >0$, there exists $\{x_1, x_2, ...x_n\} \subseteq X$ such that $A \subseteq \displaystyle\bigcup_{i=1}^{n}B(x_{i}, \epsilon)$. We call a function $f: (X,d) \rightarrow (Y,\rho)$ between two metric spaces to be totally bounded regular or TB-regular for short, if it preserves totally bounded sets, that is, for every totally bounded subset $A$ of $X$, $f(A)$ is also totally bounded. It is well-known that a metric space is totally bounded if and only if each sequence in it has a Cauchy subsequence. Intuitively, TB-regular functions can have a close relation with the functions that preserve Cauchy sequences called Cauchy continuous functions. The class of Cauchy continuous functions is well-studied \cite{[BMG],[JK2],[S2],[S1]} and lies between the class of uniformly continuous functions and that of continuous functions. TB-regular functions were first considered in \cite{[BSL]}, where the authors proved that a function is TB-regular if and only if it takes every Cauchy sequence to a sequence having a Cauchy subsequence. Thus it is immediate that every Cauchy continuous function (and hence every uniformly continuous) is TB-regular. TB-regular functions were then considered in \cite{[BL2]}, where the authors introduced an interesting class of functions that are stronger than the uniformly continuous functions, called the strongly uniformly continuous functions. TB-regular functions played a key role in the analysis of strong uniform continuity, particularly to determine the structure of the bornology of the subsets on which a (continuous) function is strongly uniformly continuous. \\
\indent The main aim of this article is to view TB-regular functions in terms of continuous functions and some well-known stronger notions of continuity. We discuss some subclasses of the class of continuous functions from a metric space $(X,d)$ to another metric space $(Y, \rho)$ such that $f: (X,d) \rightarrow (Y, \rho)$ is TB-regular if and only if whenever $f$ is followed in a composition by any function from those subclasses, it is again a TB-regular function (see Theorem \ref{thm_comp}). We note that these subclasses include the class of Cauchy continuous functions, some Lipschitz-type functions, but not the class of uniformly continuous functions. We have already seen that every Cauchy continuous (and uniformly continuous) function is TB-regular. Here we study the conditions under which the reverse implications are also true (see Theorem \ref{Cauchy_B}). In the process of studying relations between TB-regular functions and continuous functions, we establish many new useful characterizations of complete metric spaces (e.g., see Theorem \ref{Pre_2} and Theorem \ref{approx}). Finally, we give some applications of strong uniform continuity to some metric spaces possessing stronger properties than the complete metric spaces (Theorem \ref{SUC_cofinal} and Theorem \ref{SUC_UC}) and characterize the condition under which every non-vanishing strongly uniformly continuous function is stable under reciprocation (Theorem \ref{reci_suc}).
\section{\textbf{Preliminaries}}
Let us first record the precise definitions of TB-regular and Cauchy continuous functions.
\begin{definition}
	A function $f: (X,d) \rightarrow (Y,\rho)$ between two metric spaces is said to be \textit{TB-regular} if it preserves totally bounded sets, that is, for every totally bounded subset $A$ of $X$, $f(A)$ is also totally bounded.
\end{definition}
\begin{definition}
	A function $f: (X,d) \rightarrow (Y, \rho)$ between two metric spaces is said to be Cauchy-continuous if ($f(x_{n}))$ is Cauchy in $(Y,\rho)$ for every Cauchy sequence $(x_{n})$ in $(X,d)$.
\end{definition}
Throughout our analysis on TB-regular functions we use the following fact: a function is TB-regular if and only if it takes every Cauchy sequence to a sequence having a Cauchy subsequence \cite[Proposition 5.7 (1)]{[BSL]}. The following example shows that a TB-regular function need not be continuous.
\begin{example}\label{p_ex1}
	Let $X= \{\frac{1}{n}, 0: n \in \N \}$ and consider the function $f:X \rightarrow \mathbb{R}$ with $f(0) =0, ~f(\frac{1}{n}) = 1$ for odd $n$ and $f(\frac{1}{n}) = 2$ for even $n$. Clearly, $f$ is TB-regular but not even continuous. On the other hand, if we define a function $g$ on $Y = \{\frac{1}{n} : n \in \N \}$ such that $g\Big(\frac{1}{n}\Big) = n~\forall \thinspace n \in \mathbb{N}$, then $g$ is continuous but not TB-regular. 
\end{example}

Recall that a family $\mathrm{\textbf{B}}$ of non-empty subsets of a metric space $(X,d)$ is said to be \textit{bornology} if $(i)$ $\mathrm{\textbf{B}}$ forms a cover of $X$ $(ii)$ $\mathrm{\textbf{B}}$ is hereditary, that is, if $A \in \mathrm{\textbf{B}}$ and $B \subseteq A$, then $B \in \mathrm{\textbf{B}}$ $(iii)$ $\mathrm{\textbf{B}}$ is closed under finite unions (see \cite{[BSL]} and the references therein). A family of subsets $\mathcal{B}_o$ is said to be a base for $\mathrm{\textbf{B}}$ if $\forall$ $B \in \mathrm{\textbf{B}}$, $\exists \thinspace B_o \in \mathcal{B}_o$ such that $B \subseteq B_o$. If each member of a base is a closed subset of $X$, then the base is called a \textit{closed base}. As an example, note that the smallest bornology on $X$ is the family of finite subsets of $X$ and the largest one is the power set $\mathbb{P}(X)$ of $X$. Note that the family of totally bounded subsets of $X$ also forms a bornology. So in other words, a function $f: (X,d) \rightarrow (Y,\rho)$ is TB-regular if the function image of the bornology of totally bounded subsets of $(X, d)$ is contained in the bornology of totally bounded subsets of $(Y, \rho)$. 
\begin{proposition}
	Let $f: (X,d) \rightarrow (Y, \rho)$ be a function between two metric spaces. Then the set of subsets on which $f$ is TB-regular forms a bornology $\mathbf{B_T}$. Moreover, if $f$ is continuous, then $\mathbf{B_T}$ has a closed base.
\end{proposition}
\begin{proof}
	We will prove that if $f$ is TB-regular on $A$ and $B$, where $A$ and $B$ are subsets of $X$, then $f$ is TB-regular on $A \cup B$. Let $(x_n)$ be a Cauchy sequence in $A \cup B$. Thus there exists a subsequence $(x_{n_k})$ of $(x_n)$ which lies either in $A$ or in $B$. Since $f$ is TB-regular on both $A$ and $B$, there exists a Cauchy subsequence of $(f(x_n))$ in $(Y, \rho)$. Thus $f$ is TB-regular on $A \cup B$.\\
	\indent Next, let $f$ be continuous. Let $f$ be TB-regular on a subset $A$ of $X$. To prove that $f$ is TB-regular on the closure of $A$, say $\overline{A}$, let $(x_n)$ be a Cauchy sequence in $\overline{A}$. Since $f$ is continuous, for each $n \in \N$, $\exists \thinspace y_n \in A$ such that $d(x_n, y_n) < \frac{1}{n}$ and $\rho(f(x_n), f(y_n)) < \frac{1}{n}$. So $(y_n) \subseteq A$ is a Cauchy sequence. Thus $(f(y_n))$ has a Cauchy subsequence, which further implies that $(f(x_n))$ has a Cauchy subsequence too. Hence we are done.\qed
\end{proof}

In this article, we cast light on the function between metric space that are TB-regular, that is, the case when $\mathbb{P}(X) = \mathbf{B_T}$. In \cite{[BL2]}, the authors introduced an interesting class of functions that are stronger than the uniform continuous functions, called the strongly uniformly continuous functions defined as follows.
\begin{definition}
	Let $(X,d)$ and $(Y, \rho)$ be metric spaces and let $B$ be a subset of $X$. A function $f: X \rightarrow Y$ is said to be strongly uniformly continuous on $B$ if $\forall \thinspace \epsilon > 0$, there exists $\delta > 0$ such that if $d(x, y) < \delta$ and $\{x, y\} \cap B \neq \emptyset$, then $\rho(f(x), f(y)) < \epsilon$.
\end{definition}
It is interesting to note that a function $f$ on a metric space $(X,d)$ is continuous if and only if it is strongly uniformly continuous on $\{x\}$ for each $x \in X$. Also, if $f$ is a continuous function between any two metric spaces and $K$ is a compact subset of the domain, then $f$ is strongly uniformly continuous on the set $K$. In \cite{[BL2]}, TB-regular functions played a key role in the analysis of strongly uniformly continuous functions. \\
\indent Other than the above mentioned classes of continuous functions, we also study some relations of TB-regular functions with some Lipschitz-type functions. Much analysis on these functions and their applications can be found in \cite{[BGG],[BGJ],[BMG],[comp_BG],[G1],[LK1]}. Let us note their precise definitions in order of increasing size.
\begin{definition}
	For metric spaces $(X,d)$ and $(Y,\rho)$, any function $f: (X,d) \ra (Y,\rho)$ is called:
	\begin{enumerate}[label=(\alph*)]
		\item \emph{Lipschitz} if there exists $k >0$ for which $\rho(f(x), f(x')) \leq kd(x,x')~\forall \thinspace x, x' \in X$.
		\item \emph{uniformly locally Lipschitz} if there exists $\delta >0$ such that for every $x \in X$, there exists $k_x >0$ such that $\rho(f(u), f(w)) \leq k_xd(u,w)$, whenever $u, w \in B(x, \delta)$.
		\item \emph{Cauchy-Lipschitz} if $f$ is Lipschitz when restricted to the range of each Cauchy sequence $(x_n)$ in $X$.
		\item[$(e)$]\emph{locally Lipschitz} if for each $x \in X$, there exists $\delta_x >0$ such that $f$ restricted to $B(x, \delta_x)$ is Lipschitz.
	\end{enumerate}
\end{definition}
Note that the function $g$ of Example \ref{p_ex1} shows that a locally Lipschitz function need not be TB-regular, while, since every Cauchy-Lipschitz function is Cauchy continuous, every Cauchy-Lipschitz function is TB-regular.

\section{Results}
We would like to begin with the necessary and sufficient conditions under which every TB-regular function is continuous, Cauchy-continuous and uniformly continuous. The routine proof of the following result has been omitted.
\begin{theorem}
	Let $(X,d)$ be a metric space.
	\begin{enumerate}[label=(\alph*)]
		\item Every TB-regular function from $(X,d)$ to any other metric space $(Y, \rho)$ is continuous if and only if $(X,d)$ is discrete.
		\item Every TB-regular function from $(X,d)$ to any other metric space $(Y, \rho)$ is Cauchy-continuous if and only if $(X,d)$ is complete and discrete.
		\item Every TB-regular function from $(X,d)$ to any other metric space $(Y, \rho)$ is uniformly continuous if and only if $(X,d)$ is uniformly discrete.
	\end{enumerate}
\end{theorem}

 It is known that a metric space $(X,d)$ is complete if and only if each real-valued continuous function defined on it is Cauchy continuous. The next result shows that the TB-regularity of continuous functions also characterizes complete metric spaces.
\begin{theorem}\label{Pre_2}
	Let $(X,d)$ be a metric space. The following are equivalent.
	\begin{enumerate}[label=(\alph*)]
		\item $(X,d)$ is complete.
		\item Each continuous function from $(X,d)$ to any other metric space $(Y, \rho)$ is TB-regular.
		\item Each locally Lipschitz function from $(X,d)$ to any other metric space $(Y, \rho)$ is TB-regular.
		\item Each real-valued locally Lipschitz function on $(X, d)$ is TB-regular.
	\end{enumerate}
\end{theorem}
\begin{proof}  The implications $(a) \Rightarrow (b)\Rightarrow (c) \Rightarrow (d)$ are all immediate.
	
	$(d)\Rightarrow (a)$: Suppose $(X,d)$ is not complete. Thus there exists a Cauchy sequence $(x_n)$ of distinct points in $(X,d)$ such that it does not converge. Let $(\widehat{X}, \hat{d})$ be the completion of $(X,d)$ and let $(x_n)$ converges to $\hat{x} \in \widehat{X}$. Define the following function $f$ on $X$.
	$$f(x) = \frac{1}{\hat{d}(x, \hat{x})} \hspace{.5 cm}\forall \thinspace x \in X$$
	The function $f$ is not TB-regular as $\{x_n : n \in \N\}$ is totally bounded but the its image is not even bounded. To see that $f$ is locally Lipschitz, note that the function $x \mapsto \hat{d}(x, \hat{x})$ is Lipschitz and thus locally Lipschitz. It is easy to see that if a real-valued function (which is never zero) on a metric space is locally Lipschitz, then so is its reciprocal. Hence $f$ is locally Lipschitz but not TB-regular, a contradiction.\qed
\end{proof}

The following result characterizes the class of complete metric spaces in terms of the Cauchy-continuity of a subclass of continuous functions.
\begin{theorem}\label{PT1}
	Let $(X,d)$ be a metric space. The following are equivalent.
	\begin{enumerate}[label=(\alph*)]
		\item $(X,d)$ is complete.
		\item Each continuous TB-regular function from $(X,d)$ to any other metric space $(Y, \rho)$ is Cauchy continuous.
		\item Each real-valued continuous TB-regular function on $(X, d)$ is Cauchy continuous.
	\end{enumerate}
\end{theorem}
\begin{proof} Only the implication $(c)\Rightarrow (a)$ requires some justification. Suppose there exists a Cauchy sequence $(x_n)$ of distinct points in $(X,d)$ such that it does not cluster. Define the following function on the set $A =\{x_n : n \in \N\}$.
	$$f(x_n)=\left\{\begin{array}{ll}1 & \hbox{:~$n$ is odd}\\2 & \hbox{:~$n$ is even}
	\end{array}\right.$$
	By Theorem 5.1 in \cite[p. 149]{[D1]}, $f$ can be extended to a function $F:X\rightarrow(0,2)$, such that $F$ is continuous. Thus $F$ is continuous and TB-regular, but not Cauchy-continuous.\qed
\end{proof}
\begin{theorem}\label{thm_comp}
	Let $f: (X,d) \rightarrow (Y, \rho)$ be a function between two metric spaces. The following are equivalent.
	\begin{enumerate}[label=(\alph*)]
		\item $f$ is TB-regular.
		\item If $(Z, \mu)$ is a metric space and $g: (Y, \rho) \rightarrow (Z, \mu)$ is a TB-regular function, then $g \circ f$ is TB-regular.
		\item If $(Z, \mu)$ is a metric space and $g: (Y, \rho) \rightarrow (Z, \mu)$ is a Cauchy-continuous function, then $g \circ f$ is TB-regular.
		\item If $(Z, \mu)$ is a metric space and $g: (Y, \rho) \rightarrow (Z, \mu)$ is a Cauchy-Lipschitz function, then $g \circ f$ is TB-regular.
		\item If $(Z, \mu)$ is a metric space and $g: (Y, \rho) \rightarrow (Z, \mu)$ is a uniformly locally Lipschitz function, then $g \circ f$ is TB-regular.
		\item Whenever $g: (Y, \rho) \rightarrow \R$ is a uniformly locally Lipschitz function, $g \circ f$ is TB-regular.
	\end{enumerate}
\end{theorem}
\begin{proof} The implications $(a) \Rightarrow (b) \Rightarrow (c) \Rightarrow  (d) \Rightarrow (e) \Rightarrow (f)$ are all immediate.
	
	$(f) \Rightarrow (a)$: Suppose $f$ is not TB-regular. Thus there exists a Cauchy sequence $(x_n)$ in $(X,d)$ such that $(f(x_n))$ has no Cauchy subsequence. Thus there exists $\epsilon > 0$ such that by passing to a subsequence we have $\rho (f(x_n), f(x_m)) > \epsilon~\forall \thinspace n,~m \in \N.$ Define a function $g: (Y,\rho) \rightarrow \mathbb{R}$ as follows:
	$$g(y)=\left\{ \begin{array}{lll}
	n- \frac{4n}{\epsilon}d(y, f(x_n))     & : & y\in B\Big(f(x_n), \frac{\epsilon}{4}\Big) \mbox{~for some}~ n \in \N\\
	~~~~~~~	0     & :& otherwise
	\end{array}\right.$$
	The function $g$ is uniformly locally Lipschitz because $\forall \thinspace y \in Y$, $B(y, \frac{\epsilon}{4})$ intersects at most one of the balls $B(f(x_m), \frac{\epsilon}{4})$ and $g$ restricted to each ball $B(f(x_m), \frac{\epsilon}{4})$ is Lipschitz. Now for each $n \in \N$, $(g \circ f)(x_n) = n$. Thus $(x_n)$ is a Cauchy sequence such that $(g \circ f)(x_n)$  has no Cauchy subsequence, which implies that $g \circ f$ is not TB-regular, a contradiction.\qed
\end{proof}
\begin{remark}
	Note that even if every uniformly continuous function is TB-regular, the above result need not be true if we replace Cauchy continuous functions by uniformly continuous functions. As an example, let $X = \{\frac{1}{n} : n \in \N\}$. Consider the real Hilbert space $l_2$. Let $Y \subseteq l_2$ be the closed ball centered at $0$ of radius $2$. If we denote the induced metric on $Y$ by $\rho$, then $(Y, \rho)$ is a finitely chainable metric space. It is well-known that every real-valued uniformly continuous function on a metric space $(Y,\rho)$ is bounded if and only if $(Y,\rho)$ is finitely chainable \cite{[MA]}. Define the following function $f$ on $X$.
	$$f\Big(\frac{1}{n}\Big) = e_n~\forall \thinspace n \in \N$$
	where $\{e_n : n \in \N\}$ denotes the standard orthonormal basis of $l_2$. Clearly, $f$ is not TB-regular. Now, if we take any uniformly continuous function $g: (Y, \rho) \rightarrow \R$, then $g$ would be bounded, which further implies that for every such $g$, $g \circ f$ is TB-regular.
\end{remark}
\begin{theorem}\label{approx}
	Let $(X,d)$ be a metric space. The following are equivalent.
	\begin{enumerate}[label=(\alph*)]
		\item $(X,d)$ is complete.
		\item Each continuous TB-regular function from $(X,d)$ to any other Banach space $(Y, \|.\|)$ can be uniformly approximated by Cauchy-Lipschitz functions.
		\item Each real-valued continuous TB-regular function can be uniformly approximated by Cauchy-Lipschitz functions.
	\end{enumerate}
\end{theorem}
\begin{proof}
	$(a) \Rightarrow (b)$: This follows from Theorem 4.5 in \cite{[BMG]} which says that every Cauchy continuous function from $(X,d)$ to any other Banach space $(Y, \|.\|)$ can be uniformly approximated by Cauchy-Lipschitz functions.
	
	$(b) \Rightarrow (c)$: This is immediate.
	
	$(c) \Rightarrow (a)$: Suppose $(X,d)$ is not complete. By Theorem \ref{PT1}, there exists a real-valued continuous TB-regular function $f$ on $X$ which is not Cauchy-continuous. Thus there exists a Cauchy sequence $(x_n)$ in $(X,d)$ such that $(f(x_n))$ is not Cauchy. Hence for some $\epsilon >0$, by passing to a subsequence, we get $|f(x_{2k}) - f(x_{2k +1})| > \epsilon~\forall \thinspace k \in \N$. By hypothesis, there exists a Cauchy-Lipschitz function $g: (X,d) \rightarrow \R$ such that $\sup_{x \in X} |f(x) - g(x)| < \frac{\epsilon}{3}$. Since $(x_n)$ is Cauchy, there exists $M > 0$ such that $|g(x_n) - g(x_m)| \leq M d(x_n, x_m)~\forall \thinspace n, \thinspace m \in \N$. Choose $n_o \in \N$ such that $\frac{1}{n_o} < \frac{\epsilon}{3M}$. Also, there exists $k_o \in \N$ such that $\forall \thinspace k \geqslant k_o$, $d(x_{2k}, x_{2k +1}) < \frac{1}{n_o}$. Now $\forall \thinspace k \geqslant k_o$ we have
	\begin{align*}
	|f(x_{2k}) - f(x_{2k +1})| & \leq |f(x_{2k}) - g(x_{2k})| + |g(x_{2k}) - g(x_{2k+ 1})| + |g(x_{2k+ 1}) - f(x_{2k+ 1})| \\
	& < \epsilon
	\end{align*}
	We get a contradiction, thus $(X,d)$ is complete.\qed
\end{proof}
The following example shows that on a complete metric space, the set of Cauchy-Lipschitz functions need not exhaust continuous TB-regular functions on the space.
\begin{example}
	Let $X = \{\frac{1}{2^n}, 0 : n \in \N\}$. Define the following real-valued function on $X$.
	$$f(x)=\left\{\begin{array}{ll}\frac{1}{n} & :x = \frac{1}{2^n}\mbox{~for some}~ n \in \N\\0 & : x= 0
	\end{array}\right.$$
	Clearly, $f$ is continuous and TB-regular, but since $\frac{\frac{1}{n}- \frac{1}{n+1}}{\frac{1}{2^n} - \frac{1}{2^{}n+1}} = \frac{2^{n+1}}{n^2 + n}$, $f$ is not Cauchy-Lipschitz.
\end{example}

We have noticed that every Cauchy-continuous function is TB-regular but the converse is not true. The next result adds an extra condition on TB-regular functions to make it Cauchy-continuous. For proving the result, we need \textit{Efremovic Lemma} \cite[p. 92]{[B]} which says that if $(x_n)$ and $(y_n)$ are two sequences in a metric space $(X,d)$ such that $d(x_n, y_n) > \epsilon~\forall \thinspace n \in \N$, then there exists an infinite subset $\N_1$ of $\N$ such that $d(x_k, y_l) \geqslant \frac{\epsilon}{4}~\forall \thinspace k, \thinspace l \in \N_1$.
\begin{theorem}\label{Cauchy_B}
	Let $f: (X,d) \rightarrow (Y,\rho)$ be a function between two metric spaces. Then the following are equivalent.
	\begin{enumerate}[label=(\alph*)]
		\item $f$ is Cauchy-continuous.
		\item $f$ is TB-regular and the set of subsets on which $f$ is Cauchy-continuous is closed under finite unions.
		\item $f$ is TB-regular and the set of subsets on which $f$ is Cauchy-continuous forms a bornology.
	\end{enumerate}
\end{theorem}
\begin{proof}
	The implications $(a)\Rightarrow (b)\Rightarrow (c)$ are immediate.
	
	$(c)\Rightarrow (a)$: Suppose $f$ is not Cauchy-continuous. Thus there exists a Cauchy sequence $(x_n)$ in $(X,d)$ such that $(f(x_n))$ is not Cauchy. It follows that $f$ is not uniformly continuous on the set $A = \{x_n : n \in \N\}$. So for some $\epsilon >0$, there exists sequences $(y_n)$ and $(z_n)$ in $A$ such that $d(y_n, z_n) < \frac{1}{n}$ but $\rho(f(y_n), f(z_n)) > \epsilon~\forall \thinspace n \in \N$. By the Efremovic Lemma, we can assume that $\{y_n : n \in \N\} \cap \{z_n : n \in \N\} = \emptyset$. Since $(y_n) \subseteq A$, which is totally bounded, by passing to a subsequence we can assume that $(y_n)$ is Cauchy. Hence the corresponding sequence $(z_n)$ is also Cauchy. Since $f$ is TB-regular, there exists an infinite subset $\mathbb{N}_1$ of $\N$ such that $(f(y_n))_{n \in \N_1}$ is Cauchy in $Y$ and similarly, there exists an infinite subset $\mathbb{N}_2$ of $\N_1$ such that $(f(z_n))_{n \in \N_2}$ is Cauchy in $Y$. Thus $f$ is Cauchy-continuous on $\{y_n : n \in \N_2\}$ and $\{z_n : n \in \N_2\}$ but not on their union because if we enumerate the elements of $\N_2$ in the increasing order (by the well-ordering principle) say $n_1, n_2, n_3,...$, then the sequence $y_{n_1}, z_{n_1}, y_{n_2}, z_{n_2},...$ is Cauchy but its image is not. We get a contradiction. Hence $f$ is Cauchy-continuous.\qed
\end{proof}

We have the following similar characterization for the class of uniformly continuous functions.
\begin{theorem}\label{sim_UC}
	Let $f: (X,d) \rightarrow (Y,\rho)$ be a function between two metric spaces. Then the following are equivalent.
	\begin{enumerate}[label=(\alph*)]
		\item $f$ is uniformly continuous.
		\item $f$ is TB-regular and the set of subsets on which $f$ is uniformly continuous is closed under finite unions.
		\item $f$ is TB-regular and the set of subsets on which $f$ is uniformly continuous forms a bornology.
	\end{enumerate}
\end{theorem}
\begin{proof}
	Using Theorem 4.9 in \cite{[BL2]}, this can be proved in a manner similar to the proof of Theorem \ref{Cauchy_B}.\qed
\end{proof}

In view of Theorem \ref{Cauchy_B}, let us note a modified characterization of complete metric spaces.
\begin{theorem}
	Let $(X,d)$ be a metric space. Then the following are equivalent.
	\begin{enumerate}[label=(\alph*)]
		\item $(X,d)$ is complete.
		\item If $(Y,\rho)$ is a metric space and $f: (X,d) \rightarrow (Y,\rho)$ is continuous, then the set of subsets on which $f$ is Cauchy-continuous forms a bornology.
		\item If $f: (X,d) \rightarrow \R$ is continuous, then the set of subsets on which $f$ is Cauchy-continuous forms a bornology.
	\end{enumerate}
\end{theorem}
\begin{proof}
	The implication $(a)\Rightarrow (b)$ follows from the fact that every continuous function on a complete metric space is Cauchy-continuous. The implication $(b)\Rightarrow (c)$ is immediate.
	
	$(c)\Rightarrow (a)$: Suppose $(X,d)$ is not complete. Thus, in $(X,d)$, there exists a Cauchy sequence $(x_n)$ of distinct points with no cluster point. Let $A=\{x_n : n \mbox{~is odd}\}$ and $B=\{x_n : n \mbox{~is even}\}$. Since $A$ and $B$ are closed and so is $A \cup B$, by Tietze's extension theorem, let $f: (X,d) \rightarrow \R$ be a continuous function such that $f(A) = 0$ and $f(B) = 1$. Now, $f$ is Cauchy-continuous on $A$ and $B$, but not on $A \cup B$. We get a contradiction.\qed
\end{proof}

Let us note some more characterizations of complete metric spaces in terms of strong uniform continuity. Note that this characterization has originated from Proposition 4.11 in \cite{[BL2]}. 

\begin{theorem}
	Let $(X,d)$ be a metric space. Then the following are equivalent.
	\begin{enumerate}[label=(\alph*)]
		\item $(X,d)$ is complete.
		\item If $(Y,\rho)$ is a metric space and $f: (X,d) \rightarrow (Y,\rho)$ is continuous, then $f$ is strongly uniformly continuous on totally bounded subsets of $(X,d)$.
		\item If $(Y,\rho)$ is a metric space and $f: (X,d) \rightarrow (Y,\rho)$ is continuous, then $f$ is strongly uniformly continuous on each Cauchy sequence of $(X,d)$.
		\item If $f: (X,d) \rightarrow \R$ is continuous, then $f$ is strongly uniformly continuous on each Cauchy sequence of $(X,d)$.
	\end{enumerate}
\end{theorem}
\begin{proof}
	The implication $(a)\Rightarrow (b)$ follows from Proposition 4.11 in \cite{[BL2]} which says that every Cauchy-continuous function on $(X,d)$ is strongly uniformly continuous on totally bounded subsets of $(X,d)$. The implications $(b)\Rightarrow (c) \Rightarrow (d)$ are immediate.
	
	$(d)\Rightarrow (a)$: Suppose $(X,d)$ contains a Cauchy sequence $(x_n)$ of distinct points without any cluster point. By Tietze's extension theorem, let $f: (X,d) \rightarrow \R$ be a continuous function such that $f(x_n) = n~\forall \thinspace n \in \N$. Now, $(x_{2n})_{n \in \N}$ is a Cauchy sequence on which $f$ is not strongly uniformly continuous. Because otherwise, for $\epsilon = \frac{1}{2}$, there would exist $\delta > 0$ such that for each $n$, if $d(x_{2n}, y) < \delta$, then $\rho(f(x_{2n}), y) < \epsilon$, which is not possible.\qed
\end{proof}

Now we turn to metric spaces which possess stronger properties than the complete metric spaces. Cofinally complete metric spaces and UC spaces are some popular examples of such spaces.
\begin{definition}
	A sequence $(x_{n})$ in a metric space $(X,d)$ is called cofinally Cauchy if $~\forall \thinspace \epsilon > 0$, there exists an infinite subset $\mathbb{N}_\epsilon$ of $\mathbb{N}$ such that for each $n,j \in \mathbb{N}_\epsilon$, we have $d(x_n,x_j) < \epsilon.$ A metric space $(X, d)$ is said to be cofinally complete if every cofinally Cauchy sequence in $X$ clusters.
\end{definition}
It is known that a metric space $(X,d)$ is cofinally complete if and only if it is uniformly paracompact \cite{[S]}, where a metric space is said to be uniformly paracompact if for each open cover $\mathcal{V}$ of X, there exists an open refinement $\mathcal{U}$ and $\delta > 0$ such that for each $x \in X$, $B(x, \delta)$ intersects only finitely many members of $\mathcal{U}$ \cite{[R]}. In \cite{[B1]}, Beer has characterized the aforesaid spaces in the following way: A metric space $(X, d)$ is cofinally complete if and only if each sequence $(x_n)$ in $X$ satisfying $\lim\limits_{n \to \infty} \nu(x_n) = 0$ clusters, where $\nu(x)= sup\{\epsilon>0: cl(B_d(x,\epsilon))$ is compact$\}$ if $x$ has a compact neighborhood, and $\nu(x)=0$ otherwise. This geometric functional is called the \textit{local compactness functional} on $X$.\\
\indent Since cofinal completeness is a stronger form of completeness, let us study 
a family of subsets of cofinally complete metric spaces on which every continuous function must be strongly uniformly continuous. Recall that a non-empty closed subset $A$ of $X$ is said to be \textit{almost nowhere locally compact} if $\forall ~\epsilon > 0$, the set $\{a\in A: \nu(a) \geq \epsilon\}$ is compact \cite{[BM],[LK2]}. In \cite{[BM]}, it has been proved that a metric space $(X,d)$ is cofinally complete if and only if whenever $A$ is a closed subset of $X$ and $B$ is almost nowhere locally compact subset such that $A$ and $B$ are disjoint, then $d(A, B) >0$. Let us denote the set of almost nowhere locally compact subsets by $AC(X)$. 
\begin{theorem}\label{SUC_cofinal}
	Let $(X,d)$ be a metric space. Then the following are equivalent.
	\begin{enumerate}[label=(\alph*)]
		\item $(X,d)$ is cofinally complete.
		\item If $(Y,\rho)$ is a metric space and $f: (X,d) \rightarrow (Y,\rho)$ is continuous, then $f$ is strongly uniformly continuous on each member of $AC(X)$.
		\item If $(Y,\rho)$ is a metric space and $f: (X,d) \rightarrow (Y,\rho)$ is continuous, then $D_d(A,B) =0$ implies $D_\rho(f(A), f(B)) =0$, where $A \in AC(X)$ and $B$ is any non-empty subset of $X$.
	\end{enumerate}
\end{theorem}
\begin{proof}
	$(a)\Rightarrow (b)$: Since $(X,d)$ is cofinally complete, $A \in AC(X)$ implies that $A$ is compact \cite[Corollary 3.7]{[BM]} and thus every continuous function on $X$ is strongly uniformly continuous on every $A \in AC(X)$.
	
	$(b)\Rightarrow (c)$: Let $f: (X,d) \rightarrow (Y,\rho)$ be a continuous function and let $A$ and $B$ be subsets of $X$ such that $A \in AC(X)$ and $D_d(A,B) =0$. To show that $D_\rho(f(A), f(B)) =0$, let $\epsilon >0$. Since $f$ is strongly uniformly continuous on $A$, $\exists \thinspace \delta > 0$ such that if $d(x, y) < \delta$ and $\{x, y\} \cap A \neq \emptyset$, then $\rho(f(x), f(y)) < \epsilon$. Choose $a \in A$ and $b \in B$ such that $d(a,,b) < \delta$ and thus $\rho(f(a), f(b)) < \epsilon$. Consequently, $D_\rho(f(A), f(B)) =0$.
	
	$(c)\Rightarrow (a)$: Suppose $(X,d)$ is not cofinally complete. By Proposition 5.5 in \cite{[BM]}, there exists $A \in AC(X)$ and a closed subset $B$ of $X$ such that $A \cap B = \emptyset$ but $D_d(A, B) =0$. Thus for each $n \in \N$, there exists $a_n \in A$ and $b_n \in B$ such that $d(a_n, b_n) < \frac{1}{n}$. Since the sequences $(x_n)$ and $(y_n)$ cannot have any cluster point, the sets $U= \{a_n : n \in \N\}$ and $V= \{b_n : n \in \N\}$ are closed and $U$ being a subset of $A$ belongs to $AC(X)$. By Tietze's extension theorem, let $f: (X,d) \rightarrow \R$ be a continuous function such that $f(U) = 0$ and $f(V) = 1$. Now, $f$ is a continuous function such that $U \in AC(X)$ and $D_d(U, V) =0$ but $D_\rho(f(U), f(V)) \neq 0$, a contradiction. \qed
\end{proof}

Another well-known class of metric spaces which is stronger than the complete metric spaces is of UC spaces (also known as Atsuji spaces) \cite{[MA]}. 
\begin{definition}
	A metric space is called a UC space if each real-valued function defined on it is uniformly continuous.
\end{definition}
The class of UC spaces is contained in the class of cofinally complete metric spaces and it iself contains the class of compact metric spaces. Many attractive characterizations of UC spaces can be found in the survey article \cite{[KJ1]}. One can observe that $\R$ with the usual metric is cofinally complete but not UC.
\begin{theorem}\label{SUC_UC}
	Let $(X,d)$ be a metric space. Then the following are equivalent.
	\begin{enumerate}[label=(\alph*)]
		\item $(X,d)$ is UC.
		\item If $(Y,\rho)$ is a metric space and $f: (X,d) \rightarrow (Y,\rho)$ is continuous, then $f$ is strongly uniformly continuous on each subset of $X$.
		\item f $(Y,\rho)$ is a metric space and $f: (X,d) \rightarrow (Y,\rho)$ is continuous, then $f$ is strongly uniformly continuous on each closed subset of $X$.
	\end{enumerate}
\end{theorem}
\begin{proof}
	The implication $(a)\Rightarrow (b)$ follows from Theorem 5.2 in \cite{[BL2]} and the implication $(b)\Rightarrow (c)$ is immediate.
	
	$(c)\Rightarrow (a)$: It is known that a metric spaces $(X,d)$ is UC if and only if for any closed subsets $A$ and $B$ of $X$ such that $A \cap B = \emptyset$, we have $D_d(A, B) >0$ \cite{[KJ1]}. Now, rest of the proof can be done in a manner similar to the proof of $(c)\Rightarrow (a)$ in Theorem \ref{SUC_cofinal}.\qed
\end{proof}
\begin{remark}
	Using the fact that each TB-regular function sends Cauchy sequences to sequences having a Cauchy subsequence, it can be easily seen that if $f$ and $g$ are two real-valued TB-regular functions on a metric space $(X,d)$, then so are the functions $f + g$ and $f . g$, where $f. g$ denotes the pointwise multiplication of the functions $f$ and $g$.
\end{remark}
Using the above remark, we can prove the following result in a manner similar to the proof of Theorem 3.12 in \cite{[LK1]}.
\begin{theorem}
	Let $\left(X, d\right)$ be a metric space. Then the following statements are equivalent:
	\begin{enumerate}[label=(\alph*)]
		\item $\left(X, d\right)$ is complete.
		\item Whenever $f:(X,d)\rightarrow\mathbb{R}$ is a locally Lipschitz TB-regular function such that $f$ is never zero, then $\frac{1}{f}$ is also locally Lipschitz and TB-regular.
		\item Whenever $f:(X,d)\rightarrow\mathbb{R}$ is a continuous TB-regular function such that $f$ is never zero, then $\frac{1}{f}$ is also continuous and TB-regular.
	\end{enumerate}
\end{theorem}

Here we would like to mention that the stability of never-zero Lipschitz-type functions under reciprocation has been recently studied in \cite{[BGG]}. Now we study the condition under which a non-vanishing strongly uniformly continuous function on a set is stable under reciprocation, that is, if $f$ is a never-zero real-valued continuous function on $(X,d)$ such that $f$ is strongly uniformly continuous on a subset $A$ of $X$, then under what conditions $\frac{1}{f}$ would be strongly uniformly continuous on $A$? Recall that a metric space $(X, d)$ is UC if and only if each sequence $(x_n)$ in $X$ satisfying $\lim\limits_{n \to \infty} I(x_n) = 0$ clusters, where $I(x)= d(x,X\setminus \{x\})$ measures the isolation of $x$ in the space. 
\begin{definition}
	A subset $A$ of a metric space $(X,d)$ is called a UC set if every sequence $(a_n)$ in A with $\lim\limits_{n \to \infty} I(a_n) = 0$ has a cluster point in $X$.
\end{definition}
In \cite{[BL2]}, it has been proved that a subset $A$ of $X$ is UC if and only if each real-valued continuous function defined on $X$ is strongly uniformly continuous on $A$. 
\begin{theorem}\label{reci_suc}
	Let $\left(X, d\right)$ be a metric space and let $\emptyset \neq A \subseteq X$. Then the following statements are equivalent:
	\begin{enumerate}[label=(\alph*)]
		\item $A$ is a UC set.
		\item Whenever $f:(X,d)\rightarrow\mathbb{R}$ is a continuous and never zero function such that $f$ is strongly uniformly continuous on $A$, then $\frac{1}{f}$ is also strongly uniformly continuous on A.
	\end{enumerate}
\end{theorem}
\begin{proof}
	$(a)\Rightarrow (b)$: This follows from the fact that every continuous function on a UC set $A$ is strongly uniformly continuous on $A$ \cite[Theorem 5.2]{[BL2]}.
	
	$(b)\Rightarrow (a)$: Suppose $A$ is not a UC set. Thus there exists a sequence $(a_n)$ of distinct points in $A$ such that $\lim\limits_{n \to \infty} I(a_n) = 0$, but $(a_n)$ has no cluster point in $X$. Without loss of generality, let us assume that $I(a_n) < \frac{1}{n}~\forall \thinspace n \in \N$.\\
	\smallskip
	\underline{\textit{Case 1:}} The sequence $(a_n)$ is pseudo-Cauchy.\\
	Let $B =\{a_n : n \in \N\}$. Let $f$ be a function on $B$ such that $f(a_n) = \frac{1}{n}~\forall \thinspace n \in \N$. Since $f$ is a bounded uniformly continous function on $B$, it can be extended to a uniformly continuous function $\hat{f}$ on the whole space $X$ \cite{[McS]}. Define the following function.
	\begin{align*}
	g: & X \longrightarrow \R\\
	& x\longmapsto \hat{f}(x) + d(x, B)
	\end{align*}
	It is clear that $g$ is uniformly continuous and is never zero. Thus $g$ is strongly uniformly continuous on $A$, but $\frac{1}{g}$ is not uniformly continuous on $A$ as $(a_n)$ is pseudo-Cauchy and $g(a_n) = n~\forall \thinspace n \in \N$.\\
	\smallskip
	\underline{\textit{Case 2:}} The sequence $(a_n)$ is not pseudo-Cauchy.\\
	Thus there exists $\epsilon >0$ and $n_o' \in \N$ such that $d(a_n, a_m) > \epsilon~\forall \thinspace n, \thinspace m \geqslant n_o'$. Let $n_o = \max\{n_o', \frac{4}{\epsilon } + 1\}$. Since $\forall \thinspace n \in \N$, $I(a_n) < \frac{1}{n}$, $\exists\thinspace x_n \in X$ such that $x_n \neq a_n$ and $d(a_n, x_n) < \frac{1}{n}$. Let $B =\{a_n, x_n : n \geq n_o\}$. Define a function $f$ on $B$ such that $f(x_n) = \frac{1}{n}$ and $f(a_n) = \frac{1}{n^2}~\forall \thinspace n \geqslant n_o$. Again extend $f$ to uniformly continuous function $\hat{f}$ on $X$. Define the following function.
	\begin{align*}
	g: & X \longrightarrow \R\\
	& x\longmapsto \hat{f}(x) + d(x, B)
	\end{align*}
	Thus $g$ is strongly uniformly continuous on $A$ as $g$ is uniformly continuous. But $\frac{1}{g}$ is not strongly uniformly continuous on $A$ as $d(a_n, x_n) < \frac{1}{n}$ but $|g(a_n) - g(x_n)| = n(n-1)~\forall \thinspace n \geqslant n_o$.\qed
\end{proof}

\end{document}